\documentclass[12pt,twoside,a4paper]{amsart}
\usepackage[dvipsnames]{xcolor}                     
\usepackage{hyperref} 
\hypersetup{breaklinks, colorlinks,                 
    linkcolor = {Blue},                             
    citecolor = {Blue},                             
    urlcolor  = {Blue}                              %
}                                                   %
\newcommand{\urlp}[1]{\hypersetup{urlcolor = {BrickRed}}\url{#1}\hypersetup{urlcolor = {OliveGreen}}} %
\usepackage[OT2, T1]{fontenc}\usepackage[russian, english]{babel}
\usepackage[T1]{fontenc}
\usepackage[latin9]{inputenc}
\usepackage[margin=20mm,top=30mm,bottom=30mm]{geometry}
\usepackage{amsmath,amssymb,amsthm,amscd,stmaryrd,graphicx,bm,latexsym}
  
\usepackage{mathrsfs}\usepackage{mathtools}\usepackage{lmodern}\usepackage{natbib}\usepackage{enumitem}\usepackage{soul}
\newlist{parts}{enumerate}{5}       
\setlist[parts]{wide,label=\textbf{\upshape(\alph*)},ref=(\alph*),align=left,labelindent=0pt,labelsep=.5em,itemsep=.3ex,topsep=0ex,%
}
\usepackage[
final,
scrtime
]{prelim2e}
\usepackage{scalerel,stackengine}
\stackMath
\usepackage{verbatimbox} 
\usepackage{xparse}
\newlength\glyphwidth
\newlength\widthofx
\newsavebox\hatglyphCONTENT
\sbox\hatglyphCONTENT{%
    \addvbuffer[-0.05ex -1.3ex]{$\hat{\phantom{.}}$}%
}
\newcommand\hatglyph{\resizebox{0.6\widthofx}{!}{\usebox{\hatglyphCONTENT}}}
\newcommand\shifthat[2]{%
    \stackengine{0.2\widthofx}{%
        \SavedStyle#2}{%
        \rule{#1}{0ex}\hatglyph}{O}{c}{F}{T}{S}%
}
\ExplSyntaxOn
\newcommand\relativeGlyphOffset[1]{%
    \str_case:nnF{#1}{%
        {A}{0.18}%
        {B}{0.1}%
        {W}{0.02}%
        {J}{0.18}%
        {\phi}{0.17}%
    }{0.05}
}\ExplSyntaxOff
\NewDocumentCommand{\hatt}{mO{#1}}{%
    \ThisStyle{%
        \setlength\glyphwidth{\widthof{$\SavedStyle{}\longleftarrow$}}%
        \setlength\widthofx{\widthof{$\SavedStyle{}x$}}%
        \shifthat{\relativeGlyphOffset{#2}\glyphwidth}{#1}%
  }%
}
\renewcommand{\PrelimText}{\footnotesize[\,
\texttt{\jobname.tex}\hfill \today\ at \thistime\,]}
  \theoremstyle{plain}
    \newtheorem{Thm}{Theorem}[section]  
      \newtheorem{Lem}[Thm]{Lemma}
   
    \theoremstyle{definition} 

\renewcommand{\cite}{\citet}
\newcommand*{\nameadjunct}{\relax}
\makeatletter
\renewcommand*{\NAT@nmfmt}[1]{\NAT@up #1\nameadjunct}
\makeatother

\newcommand*{\citeposs}[2][]{%
    \begingroup
    \renewcommand*{\nameadjunct}{'s}%
    \citet[#1]{#2}%
    \endgroup
}

\newcounter{aequiv}

\newenvironment{Aequiv}
  {\begin{list}{\rm(\roman{aequiv})}%
   {\usecounter{aequiv}%
     \setlength{\leftmargin}{0pt}%
     \setlength{\labelwidth}{2.8em}%
     \setlength{\itemindent}{1.4\labelwidth}%
     \setlength{\partopsep}{0pt}%
     \setlength{\topsep}{0pt}%
     \setlength{\labelsep}{1.3em}%
     \setlength{\itemsep}{0.0ex}\setlength{\parsep}{0.0ex}%
   }}%
 {\end{list}} 
\let\oldsqrt\sqrt
\def\sqrt{\mathpalette\DHLhksqrt}
\def\DHLhksqrt#1#2{%
\setbox0=\hbox{$#1\oldsqrt{#2\,}$}\dimen0=\ht0
\advance\dimen0-0.2\ht0
\setbox2=\hbox{\vrule height\ht0 depth -\dimen0}%
{\box0\lower0.4pt\box2}}

\newcommand{\N}{{\mathbb N}}\newcommand{\R}{{\mathbb R}}%
\renewcommand{\epsilon}{\varepsilon}\renewcommand{\phi}{\varphi} %
\renewcommand{\rho}{\varrho}\renewcommand{\theta}{\vartheta}     %
\DeclareMathOperator*{\bigconv}{\mbox{\LARGE$\ast$}}
\DeclareFontEncoding{LS2}{}{\noaccents@}                                                 %
\DeclareFontSubstitution{LS2}{stix}{m}{n}                                                %
\DeclareSymbolFont{largesymbols_stix}{LS2}{stixex}{m}{n}                                 %
\DeclareMathDelimiter{\lpp}{\mathopen}{largesymbols_stix}{"DE}{largesymbols_stix}{"02}   %
\DeclareMathDelimiter{\rpp}{\mathclose}{largesymbols_stix}{"DF}{largesymbols_stix}{"03}  %
\renewcommand{\[}{\begin{eqnarray*}}\renewcommand{\]}{\end{eqnarray*}}
\newcommand{\la}{\begin{eqnarray}}\newcommand{\al}{\end{eqnarray}}
  
\renewcommand{\Prob}{\mbox{\rm Prob}}      
\newcommand{\dd}{\mathrm{d}}\newcommand{\ee}{\mathrm{e}}\newcommand{\ii}{\mathrm{i}}

\begin{document}
\author[Mattner]{Lutz Mattner}
\address{Universit\"at Trier, Fachbereich IV -- Mathematik, 54286~Trier, Germany}
\email{mattner@uni-trier.de}
\title[]{Extreme expectations of Bernoulli convolutions\\
 given their first few  moments
 are attained\\ at shifted convolutions of as few  binomials}
\begin{abstract} 
 A result of Chebyshev~(\citeyear{Tchebichef1864}) and \citet{Hoeffding1956}, 
 on bounding an expectation of a given function with respect to 
 a Bernoulli convolution (also called Poisson binomial law, 
 or law of the number of successes in independent trials) with any given first moment, 
 is  here generalised to the case of any given first few  moments, 
 as indicated in the title.
 
 A nonprobabilistic, and perhaps more obvious, reformulation is: Every permutation invariant 
 and separately affine-linear function of $n$ real variables $x_i\in[a,b]$ assumes its extremal 
 values given the power sums $\sum_{i=1}^nx_i^1,\ldots, \sum_{i=1}^nx_i^r$ 
 at vectors $x$ with at most $r$ coordinate values different from $a$ and $b$.
\end{abstract} 
\date{\today} 
\maketitle

\section{Introduction and result}
{\em Bernoulli convolutions}, by which we mean the laws 
\la                                     \label{Eq:Def_BC}
 \mathrm{BC}_p &\coloneqq& \bigconv_{i\in I}\mathrm{B}_{p_i} 
   \quad\text{ for }p\in[0,1]^I\text{ with $I$ finite}
\al
and with each $\mathrm{B}_{p_i}$ being Bernoulli 
(
details about our terminology and notation  
are mostly given below 
between \eqref{Eq:Def_underline_n} and Theorem~\ref{Thm:Chebyshev-Hoeffding_k}),
occur naturally for two quite different reasons:

First, $\mathrm{BC}_p$ as in~\eqref{Eq:Def_BC} is the law of the number of successes 
in $\# I$ independent trials with possibly varying success probabilities $p_i$.
As such it has been studied by 
\citet[\S38, pp.~421--423, the case of $\nu^{(i)}=1$;
\citeyear{Laplace1886}, pp.~430--432]{Laplace1812},
who derived a normal approximation based on the mean and the variance 
\la                                      \label{Eq:mean_var_BC}
 \mu(\mathrm{BC}_p) \,\ =\,\ \sum_{i\in I} p_i \,,
 && \sigma^2(\mathrm{BC}_p ) \,\ =\,\ \sum_{i\in I}p_i(1-p_i)\,,
\al
and again by  \citet{Poisson1837} as described by \citet[in particular pp.~570, 579]{Hald1998}.
Because of Poisson's work, 
$\mathrm{BC}_p$ is often called a {\em Poisson binomial} law, 
to distinguish it from the special case of an ordinary binomial law 
\la                                \label{Eq:Def_binomial}
 \mathrm{B}_{n,p} &\coloneqq& \mathrm{B}_p^{\ast n} \,\ \coloneqq\,\ 
  \bigconv_{i=1}^n\mathrm{B}_{p} \,\ =\,\ \mathrm{BC}_{(p,\ldots,p)} 
  \quad\text{ for }n\in\N_0\text{ and }p\in[0,1].
\al

Second, many combinatorial probability distributions unexpectedly turn out
to be Bernoulli convolutions, with in general the $p_i$ in~\eqref{Eq:Def_BC} 
neither explicitly  known nor 
rational, 
but nevertheless the mean and the variance, and perhaps also higher moments or cumulants,
easy to obtain. 
For example, each hypergeometric law $\mathrm{H}_{n,r,b}$ of the number of red balls obtained
by drawing $n$~balls successively without replacement from an urn containing $r$ red and $b$ blue balls, 
for which \citet[pp.~729--733]{MattnerSchulz2018}
recall some basic facts and references, 
is a Bernoulli convolution, as proved by \citet[Corollary 5 with $n=2$]{VatutinMikhailov1983}.
Here it is easy to check that in particular 
\[
 \mathrm{H}_{2,r,b} &=& \mathrm{B}_{\frac{r}{r+b}+\epsilon} \ast \mathrm{B}_{\frac{r}{r+b}-\epsilon}
 \ \text{ with }\ \epsilon\coloneqq \tfrac{1}{r+b}\sqrt{\tfrac{rb}{r+b-1}}
 \quad\text{ for }r,b\in\N, 
\]
with the Bernoulli success parameters irrational, and then without any combinatorial meaning, 
for example if $r=b=2$. 
\citet[p.~279]{Pitman1997} gives references to further combinatorial Bernoulli convolutions, 
there named {\em P\'olya frequency distributions}. 
\citet{Warren1999} treats recursion techniques for proving Poisson-binomiality, 
applicable in particular to give another proof for the hypergeometric
case mentioned above. 
\citet{TangTang2019} provide a more recent review of Bernoulli convolutions. 

Hence, in either case, bounds for probabilities $\mathrm{BC}_p(A)$ 
or more general expectations $\mathrm{BC}_p\, g$, 
with a given set $A$ or function $g$ and given any first few moments or cumulants 
of $\mathrm{BC}_p$, are of interest. 
Theorem~\ref{Thm:Chebyshev-Hoeffding_k}, stated below after explaing some terminology and notation
and recalling some standard facts, 
reduces the search for such bounds to certain special Bernoulli convolutions, namely 
shifted convolutions of few 
binomial laws as in~\eqref{Eq:shifted_convolution_of_k_binomials}.

We put                            
$\N\coloneqq\{1,2,\ldots\}$, $\N_0\coloneqq\N\cup\{0\}$, and
\la                                            \label{Eq:Def_underline_n}
 \underline{n}\ \coloneqq\ \{1,\ldots,n\} \quad \text{ and } \quad  
 \underline{n}_0\ \coloneqq\ \underline{n}\cup\{0\}
  \ =\ \{0,\ldots,n\} \quad \text{ for }n\in\N_0\,, 
\al
so in particular $\underline{0}=\emptyset$.
We write $\#A$ for the number of elements of a set $A$, 
and $]a,b[$ for an open interval.

We write $\Prob(\R)$ for the set of all probability measures (or {\em distributions}, or {\em laws})
on the Borel sets of $\R$, and $\ast$ for convolution  on $\Prob(\R)$.
The special laws occuring in this paper are Dirac laws $\delta_a$ for $a\in\R$, 
Bernoulli laws $\mathrm{B}_p \coloneqq (1-p)\delta_0+p\delta_1$ for $p\in[0,1]$,
and the convolutions~$\mathrm{BC}_p$ and $\mathrm{B}_{n,p}$ 
defined by~\eqref{Eq:Def_BC} and~\eqref{Eq:Def_binomial}, 
where of course $\bigconv_{i\in\emptyset} P_i\coloneqq \delta_0$. 

For $r\in\N_0$ we put $\Prob_r(\R)\coloneqq \{P\in\Prob(\R): \int |x|^r\,\dd P(x)<\infty\}$, 
and $\mu_r(P)\coloneqq\int |x|^r\,\dd P(x)$ for $P\in\Prob_r(\R)$;
then $\mu\coloneqq \mu_1$ on $\Prob_1(\R)$, $\sigma^2\coloneqq \mu_2-\mu_1^2$ on $\Prob_2(\R)$.  
For $r\in\N$, we let $\kappa_r$ denote the $r$th cumulant functional on $\Prob_r(\R)$, 
that is, 
\[
 \kappa_r(P) &\coloneqq& \ii^{-r}  
 \Big(\frac{\dd}{\dd t}\Big)^{\!r} \!
  \log\big(\int\ee^{\ii tx}\,\dd P(x)\big)\Big|_{t=0}
  \quad\text{ for } P \in \Prob_r(\R).
\]
We recall \citeposs{Thiele1889,Thiele1899} recursion
\la                                            \label{Eq:Thiele_recursion}
 \mu^{}_{r+1} &=& \sum_{\ell=0}^r\textstyle{\binom{r}{\ell}} \mu^{}_{r-\ell}\kappa^{}_{\ell+1}
 \quad\text{ on } \Prob_{r+1}(\R), 
 \text{ for each }r\in\N_0
\al
given in \citet[pp.~144, 150]{Hald2000} with essentially its modern proof as, 
for example, in \citet{Mattner1999};
this shows that the first $r$ moments $\mu_1,\ldots, \mu_r$ 
are (polynomial) functions of the first $r$ cumulants $\kappa_1,\ldots, \kappa_r$
and conversely;
in particular $\kappa_1=\mu$, $\kappa_2=\sigma^2$. Hence it is purely a matter of taste whether 
one fixes some first few moments or the same number of first cumulants 
of a law on $\R$, and we have chosen the latter in 
Theorem~\ref{Thm:Chebyshev-Hoeffding_k} below, and the perhaps 
better known former in the title of this paper.
The main argument in favour of cumulants is of course their additivity with respect 
to convolution, 
\la                                                           \label{Eq:cumulants_additive}
    \kappa_r(P\ast Q)&=& \kappa_r(P)+\kappa_r(Q) \quad\text{ for }r\in\N\text{ and }P,Q\in\Prob_r(\R), 
\al
first observed by \citet{Thiele1889} according to  \citet[pp.~345, 347]{Hald1998}, 
and later shown to be characteristic in a quite strong sense by
\citet{Mattner1999,Mattner2004}.

We use 
the notation $Pg\coloneqq \int g\,\dd P$ for $P$-integrable functions $g$.

\begin{Thm}                       \label{Thm:Chebyshev-Hoeffding_k}
Let $I$ be a set with $n\coloneqq \#I\in\N_0$, $g: \underline{n}_0
\rightarrow \R$ a function, $r\in\N_0$, 
and $c \in \R^r$ with 
\la            \label{Eq:Def_D_nonempty}
 D&\coloneqq& \big\{p\in[0,1]^I : \kappa_\ell(\mathrm{BC}_p)=c_\ell 
  \text{ for }\ell\in \underline{r} \big\} 
  \,\ \neq\,\ \emptyset \,.
\al
Then $\max\limits_{p\in D} \mathrm{BC}_p\, g$ 
is attained at some $p\in D$
with $r' \coloneqq\#(\{p_i:i\in I \}\!\setminus\!\{0,1\}) \le r$.
For such a~$p$, we have
\la                           \label{Eq:shifted_convolution_of_k_binomials}
 \mathrm{BC}_p &=& \delta_{n_0} \ast \bigconv_{j=1}^{r'} \mathrm{B}_{n_j,q_j}
\al
with $\{q_j : j\in\underline{r'}\}\coloneqq\{p_i:i\in I \}\!\setminus\!\{0,1\} $, 
$q_0\coloneqq 1$, and $n_j \coloneqq \#(\{i\in I: p_i=q_j\})$ for $j\in
\underline{r'}_0$.
\end{Thm}

The corresponding minimisation result is of course obtained by applying the above to $-g$.

For $r=0$, Theorem~\ref{Thm:Chebyshev-Hoeffding_k} is a rather obvious and not very interesting 
consequence of the separate affine-linearity of 
$\mathrm{BC}_p\, g$ as a function of $p\in D=[0,1]^I$.
For $r=1$, Theorem~\ref{Thm:Chebyshev-Hoeffding_k} 
is \citet[p.~717, Corollary~2.1]{Hoeffding1956}, 
proved explicitly before by \citet[p.~260, second Th\'eor\`eme]{Tchebichef1864} 
for $g$ an indicator of a set of the form $\{x\in\N_0 : x\ge m\}$, 
but with a method applicable to general $g$ as well;
examples of its application, with the determination of an appropriate pair $(n_0,n_1)$ not quite trivial,
include \citet[p.~202]{MattnerRoos2007} and \citet[p.~307]{MattnerTasto2015}.

We hope that the present Theorem~\ref{Thm:Chebyshev-Hoeffding_k}
will find comparable applications for $r\ge2$, perhaps in the following nonprobabilistic reformulation.
For $r\in\N_0$, $I$ a finite set, and $x\in\R^I$, we denote 
the $r$th power sum of $x$ by 
\la                                             \label{Eq:Def_power_sum}
 S_r(x) &\coloneqq& \sum_{i\in I}x_i^r \,. 
\al

\begin{Thm}                          \label{Thm:Extreme_symmetric_multiaffine_function_given_per_sums}
Let $I$ be a finite set, $a,b\in\R$, $f:[a,b]^I\rightarrow \R$ a function permutation invariant
in its $\#I$ variables and  affine-linear in each of these, $r\in\N_0$, and 
$c\in\R^r$ with 
\[
 D &\coloneqq& \{ x\in[a,b]^I : S_\ell(x) = c_\ell\text{ for }\ell\in\underline{r}\}
  \,\ \neq \,\ \emptyset.
\]
Then $\max_{x\in D} f(x)$ is attained at some $x\in D$ with $\#(\{x_i:i\in I\}\!\setminus\!\{a,b\})\le r$.
\end{Thm}

Theorem~\ref{Thm:Extreme_symmetric_multiaffine_function_given_per_sums} is closely 
related to \citet[p.~221, Theorem~3, 
there given without details of a proof]{KovacecKuhlmannRiener2012}.
In fact, the present Lemma~\ref{Lem:S_{r+1}_extremal} is, 
after the first two sentences of its proof, a special case of the claim just cited.

The proof of Theorem~\ref{Thm:Chebyshev-Hoeffding_k}
given near the end of the next section rests on 
Lemmas~\ref{Lem:S_{r+1}_extremal}--\ref{Lem:Symmetric_partially_affine_functions}.
These are ``well--known'', 
as just indicated in case of Lemma~\ref{Lem:S_{r+1}_extremal}, 
but we attempt to 
provide appropriate references and, if deemed advisable, details of proofs.
The idea for Theorem~\ref{Thm:Chebyshev-Hoeffding_k} is then to apply Lemma~\ref{Lem:S_{r+1}_extremal} to  
any $r+1$ coordinates of a suitable maximiser~$p$. 

\section{Proofs}
\begin{Lem}                            \label{Lem:S_{r+1}_extremal}
Let $I$ be a finite set, $a,b\in\R$, $r\in\N_0$, $c\in\R^r$,  and 
let $x$ locally minimise or maximise  
$S_{r+1}$ in $[a,b]^I$ subject to the constraints $ S_\ell(x)=c_\ell$ 
for $\ell\in\underline{r}$. Then
\la                           \label{Eq:At_most_k_nonboundary_x_i}
   \#(\{x_i:i\in I\}\!\setminus\!\{a,b\}) &\le& r\,.
\al
\end{Lem}
\begin{proof}
Let $I_0\coloneqq \{i \in I : x_i\in\mathopen]a,b\mathclose[ \}$.
Then $(x_i:i\in I_0)$ locally minimises or maximises $\sum_{i\in I_0} x_i^{r+1}$ 
in the open set $]a,b[^{I_0}$ subject to 
the constraints $\sum_{i\in I_0}x_i^\ell=d_\ell\coloneqq c_\ell- \sum_{i\in I\setminus I_0}x_i^\ell$ 
for $\ell\in\underline{r}$. 
Hence there exists a Lagrange multiplier 
$\lambda \in \R^{\underline{r}_0}_{}$, not identically zero, with
\[
 \lambda_0(r+1)x_i^r+ \sum_{\ell=1}^r \lambda_\ell \ell x_i^{\ell-1} &=& 0
  \quad \text{ for }i\in I_0\,,
\]
for example by \citet[p.~63]{Tichomirov1982}.
Hence the $x_i$ with $i\in I_0$ are zeros of a common nonzero polynomial of degree at most $r$,
hence $\#\{x_i:i\in I_0\}\le r$, and hence we have~\eqref{Eq:At_most_k_nonboundary_x_i}. 
\end{proof}

\begin{Lem}                            \label{Lem:Cumulants_of_B_p}
For $r\in\N$, there are unique $a_{r,\ell}\in\R$ with 
\la                                                     \label{Eq:Cumulants_of_B_p}
 \kappa_r(\mathrm{B}_p) &=& \sum_{\ell=1}^r a_{r,\ell}\,p^\ell \quad\text{ for }p\in[0,1],
\al
and then $a_{r,r}=(-1)^{r-1}(r-1)!$.
In particular we have $\kappa_1(\mathrm{B}_p)= p  $, $\kappa_2(\mathrm{B}_p)= p(1-p)  $, 
$\kappa_3(\mathrm{B}_p)=p(1-p)(1-2p)   $, $\kappa_4(\mathrm{B}_p)= p(1-p)\big(1-6 p(1-p)\big)  $.
\end{Lem}
\begin{proof}[First proof]
Since obviously $\mu_r(\mathrm{B}_p)=p^{1\wedge r}$ for $r\in\N_0$ and $p\in[0,1]$,
applying~\eqref{Eq:Thiele_recursion} to $\mathrm{B}_p$ yields 
\[
 \kappa_{r+1}(\mathrm{B}_p) 
   &=& p\,\Big(1- \sum_{\ell=0}^{r-1}\textstyle{\binom{r}{\ell}} \kappa_{\ell+1}(\mathrm{B}_p)\Big)
   \quad\text{ for }r\in\N_0\text{ and }p\in[0,1],
\]
and hence the claim by induction.
\end{proof}
\begin{proof}[Second proof]
One easily verifies
\[                                       
 \kappa_1(\mathrm{B}_p)\ =\ p,
 && \kappa_{r+1}(\mathrm{B}_p) \ = \ p(1-p)\frac{\dd}{\dd p}\kappa_r(\mathrm{B}_p)
 \quad \text{ for }r\in\N\text{ and }p\in[0,1]
\]
as indicated in \citet[p.~203, Exercise 5.1]{Kendall}, 
with reference to \citet{Frisch1926}
and \citet[pp.~392--393, there $p$ and $1\!-\!p$ inconsistently interchanged]{Haldane1939or1940}.
\end{proof}

We fix the notation $a_{r,\ell}$\,,  defined by~\eqref{Eq:Cumulants_of_B_p},
up to the proof of Lemma~\ref{Lem:S,K,E}.
In generalisation of~\eqref{Eq:mean_var_BC} we hence get:

\begin{Lem}
For $r\in\N$ we have 
\la                          \label{Eq:Cumulants_of_BC_p}
 \kappa_r(\mathrm{BC}_p) &=& \sum_{i\in I} \kappa_r(\mathrm{B}_{p^{}_i})
   \,\ =\,\ \sum_{\ell=1}^r a_{r,\ell} S_\ell(p) \quad\text{ for }p\in[0,1]^I
 \text{ with $\#I$ finite}. 
\al
\end{Lem}
\begin{proof}
First 
\eqref{Eq:cumulants_additive},
then~\eqref{Eq:Cumulants_of_B_p}, then~\eqref{Eq:Def_power_sum}.
\end{proof}

For $r\in\N_0$, $I$ a finite set, and $x\in\R^I$, we put
\la
 K_r(x) &\coloneqq& \sum_{\ell=1}^r a_{r,\ell} S_\ell(x) \,, \label{Eq:Def_K_r} \\
 E_r(x) &\coloneqq& \sum_{\alpha\in\{0,1\}^I : |\alpha|=r} x^\alpha     \label{Eq:Def_E_r}
   \,\ =\,\ \sum_{J\subseteq I: \#J=r}\ \prod_{i\in J}x_i \,, 
\al 
using here multi-index notation, as explained for example by \citet[pp.~54--56]{John1982}.
So the $K_r$ are, in case of $r\ge 1$, the polynomial extensions of the left hand sides
of~\eqref{Eq:Cumulants_of_BC_p}, and the $E_r$ are the elementary symmetric functions, 
$E_0(x)=1$, $E_1(x)=\sum_{i\in I}x_i=S_1(x)$,  \ldots, $E_{\#I}(x)=\prod_{i\in I}x_i$,
$E_r(x)=0$ for $r>\#I$.
 
Referring to the use of the symbols $S,K,E$ 
in~(\ref{Eq:Def_power_sum},\ref{Eq:Def_K_r},\ref{Eq:Def_E_r}), 
we have:

\begin{Lem}                       \label{Lem:S,K,E}
Let $r\in\N$. For $(A,B) \in\{S,K,E\}^2$ there exist 
a constant $c^{}_{A,B,r}\in\R\!\setminus\!\{0\}$ and a polynomial function 
$T_{A,B,r}:\R^{r-1} \rightarrow \R$ with 
\[
  A_r(x) &=& c^{}_{A,B,r} B_r(x) + T_{A,B,r}(B_1(x),\ldots, B_{r-1}(x))
  \quad \text{ for } x\in\R^I\text{ with }I\text{ finite}.
\]
In particular, $c^{}_{E,S,r}=\frac{(-1)^{r-1}}{r}$.
\end{Lem}
%
%
\begin{proof}
For $(A,B)=(K,S)$, the claim is true with $c^{}_{K,S,r}=a_{r,r}= (-1)^{r-1}(r-1)!$
and $T_{K,S,r}$ even linear, by~\eqref{Eq:Def_K_r} and Lemma~\ref{Lem:Cumulants_of_B_p}.

For $(A,B)=(E,S)$, the claim follows from the Newton identities 
\[
 E_r &=& \frac{1}{r} \sum_{\ell=1}^{r} (-1)^{\ell-1} S_{\ell} E_{r-\ell} \quad\text{ for }r\in \N
\]
proved, for example, by \citet[p.~261]{Uspensky1948} and by \citet[p.~23, (2.11$'$)]{Macdonald1995}.

The other 
cases follow from these two, with $c^{}_{A,C,r}=c^{}_{A,B,r}c^{}_{B,C,r}$
for $A,B,C\in\{S,K,E\}$. 
\end{proof}

As noted by \citet[pp.~713-714]{Hoeffding1956}, we have:
\begin{Lem}                            \label{Lem:Symmetric_partially_affine_functions}
Let $I$ be a finite set with $n\coloneqq \#I\in\N_0$, and let  
and $f:[0,1]^I\rightarrow \R$ be a function.
Then the following three statements are equivalent:
\begin{Aequiv}
\item There is a function $g:\underline{n}_0\rightarrow\R$  \label{Aequiv:f_BC_expectation} 
 with $f(p)=\mathrm{BC}_p\, g$ for $p\in[0,1]^I$.  
\item $f$ is permutation invariant in its $n$     \label{Aequiv:f_partially_affine_etc}
 variables, and affine-linear in each of these. 
\item There is a vector $b\in \R^{\underline{n}_0}_{}$ with   \label{Aequiv:f_lin_comb_of_E_k}
 $f(p) = \sum_{k=0}^{n} b_k E_k(p)$  for $p\in[0,1]^I$.
\end{Aequiv}   
\end{Lem}
\begin{proof}
 \ref{Aequiv:f_BC_expectation}  $\Rightarrow$ \ref{Aequiv:f_partially_affine_etc}: Obvious
 by permutation invariance and multilinearity of convolution. 
 
 \smallskip \ref{Aequiv:f_partially_affine_etc} $\Rightarrow$ \ref{Aequiv:f_lin_comb_of_E_k}:
 If 
 $f$ is separately affine-linear as stated, then $f$ is in particular polynomial 
 (jointly in its variables, not merely separately), as can be seen by induction:
 Assuming here $I=\underline{n}$ and $n\ge 1$, we have
 $f(p)=f(p_1,\ldots,p_{n-1},0) + (f(p)- f(p_1,\ldots,p_{n-1},0))
  =  g(p_1,\ldots,p_{n-1})+ h(p_1,\ldots,p_{n-1})p_n$ for $p\in[0,1]^n$ 
  with some functions $g$ and $h$, which by considering for example first $p_n=0$ and then $p_n=1$, 
 are seen to be affine-linear in each of their respectively $n-1$ variables, 
 hence polynomial by assumption, yielding polynomiality of $f$.
 
 Hence then $f$ is representable by its Taylor expansion about $0\in[0,1]^I$, 
 which under permutation invariance and separate affine-linearity reduces to the form claimed
 in~\ref{Aequiv:f_lin_comb_of_E_k}.
 
 \smallskip \ref{Aequiv:f_lin_comb_of_E_k} $\Rightarrow$ \ref{Aequiv:f_BC_expectation}:
 For $k\in\underline{n}_0$ and $p\in[0,1]^I$, we get 
 \[
  E_k(p) &=& \int\limits_{\{0,1\}^I}E_k\,\dd\bigotimes_{i\in I}\mathrm{B}_{p^{}_i}
   \,\ =\,\ \int\limits_{\{0,1\}^I}\binom{S_1}{k}\,\dd\bigotimes_{i\in I}\mathrm{B}_{p^{}_i}
   \,\ =\,\ \int\limits_{\underline{n}_0} \binom{x}{k}\dd \mathrm{BC}_p(x)
 \]
 by using in the third step that the convolution $\mathrm{BC}_p$ is by definition
 the image under $S_1$ of the product measure $\bigotimes_{i\in I}\mathrm{B}_{p^{}_i}$, 
 in the second $E_k(y)=\binom{S_1(y)}{k}$ for $y\in\{0,1\}^I$, 
 and in the first $\int\sum=\sum\int$ and Fubini. 
 Hence, given  \ref{Aequiv:f_lin_comb_of_E_k}, we take $g(x)\coloneqq\sum_{k=0}^{n} b_k \binom{x}{k}$
 for $x\in\underline{n}_0$.
\end{proof}

Instead of the  proof of polynomiality
in \ref{Aequiv:f_partially_affine_etc} $\Rightarrow$ \ref{Aequiv:f_lin_comb_of_E_k}
above, it would appear to be more elegant to use a general fact 
like the generalisation by \citet[Theorem in section 5, Theorem 8.3]{Palais1978}
of a result of \citet{Carroll1961}, but we are not aware of a short exposition yielding easily 
what is needed here.

\begin{proof}[Proof of Theorem~\ref{Thm:Chebyshev-Hoeffding_k}]
Under the assumptions of the theorem up to~\eqref{Eq:Def_D_nonempty}, 
$\max_{p\in D} \mathrm{BC}_p\, g$ exists by continuity, compactness, and nonemptyness.
So we may choose a maximiser of~$S_{r+1}$ 
among the perhaps several maximisers of $D \ni p\mapsto \mathrm{BC}_p\, g$
and, changing notation, call this maximiser~$\hatt{p}$, rather 
than $p$ as in the theorem.

We recall that an $x\in\R^I$ is a function on $I$, hence 
a set $\{(i,x_i):i\in I\}$ of ordered pairs;
hence the notation $\hatt{p}|_J\coloneqq\{(j,\hatt{p}_j): j\in J\}$ 
and $q\cup \hatt{s}$ occurring below.

Let us assume that $J\subseteq I$ with $\# J=r+1$, 
and let $\hatt{q}\coloneqq \hatt{p}|_J$ and    
$\hatt{s}\coloneqq \hatt{p}|_{I\setminus J}$. 
With 
\la                                                   \label{Eq:Def_f}
 f(q) \coloneqq \mathrm{BC}_{q\cup \hatt{s}}\,g  \quad\text{ for } q\in [0,1]^J 
\al
then $\hatt{q}$ maximises $f$ on 
\[
 D_{\hatt{s}} &\coloneqq&  \big\{ q\in[0,1]^J : q\cup \hatt{s} \in D\big\}
  \,\ =\,\  \big\{q \in[0,1]^J : \kappa_\ell(\mathrm{BC_q}) 
    = c_\ell - \kappa_\ell(\mathrm{BC}_{\hatt{s}})    \text{ for }\ell\in\underline{r}\big\}\\
  &=&   \big\{q \in[0,1]^J : S_\ell(q)=d_\ell  \text{ for }\ell\in\underline{r}\big\}
\]
with a certain $d=d(c,\hatt{s})\in\R^r$. Here 
the second representation of $D_{\hatt{s}}$  follows from 
$\mathrm{BC}_{q\cup \hatt{s}} =   \mathrm{BC}_{q}\ast\mathrm{BC}_{\hatt{s}}$
and~\eqref{Eq:cumulants_additive}, 
and the third from (\ref{Eq:Cumulants_of_BC_p},\ref{Eq:Def_K_r}) 
and Lemma~\ref{Lem:S,K,E} applied to $(S,K)$ and to $(K,S)$, 
with $q$ in the role of $x$, and each $\ell\in\underline{r}$ here in the 
role of $r$ there.

Now $f$ is a permutation invariant function of its $r+1$ arguments, 
and affine-linear in each of these, and hence 
the implication \ref{Aequiv:f_partially_affine_etc} $\Rightarrow$ \ref{Aequiv:f_lin_comb_of_E_k}
in Lemma~\ref{Lem:Symmetric_partially_affine_functions}, followed by
Lemma~\ref{Lem:S,K,E} applied to $(E,S)$ and  to each  $\ell\in\underline{r+1}$ here in the 
role of $r$ there, yields
\la  \label{Eq:f_via_S_1_..._S_{r+1}}
 \quad f(q) &=& \sum_{\ell=0}^{r+1} b_\ell E_\ell(q) 
  \,\ =\,\ \frac{(-1)^{r}b_{r+1} }{r+1}S_{r+1}(q) + h(S_1(q),\ldots,S_r(q))
 \quad\text{ for } q\in[0,1]^J
\al
for some $b=b(\hatt{s})\in \R^{\underline{r+1}_0}_{}$ and some function $h=h_b$. 

If $\hatt{q}$ is a maximiser or minimiser of  $S_{r+1}$ on $D_{\hatt{s}}$, 
then $\#\{i\in J: \hatt{q}_i\in\mathopen]0,1\mathclose[\,\}\le r$ by Lemma~\ref{Lem:S_{r+1}_extremal}.
If $\hatt{q}$ is neither, then~\eqref{Eq:f_via_S_1_..._S_{r+1}} 
and $f(\hatt{q})=\max_{q\in D_{\hatt{s}}}f(q)$ imply $b_{r+1}=0$, 
hence $f(q)=f(\hatt{q})$ for every $q\in D_{\hatt{s}}$, 
so that by~\eqref{Eq:Def_f} each $q\cup\hatt{s}$ with $q\in D_{\hatt{s}}$ maximises $D \ni p\mapsto \mathrm{BC}_p\, g$, 
and hence  the second maximality property of $\hatt{p}$ yields
%
$S_{r+1}(\hatt{q})=S_{r+1}(\hatt{p})-S_{r+1}(\hatt{s}) \ge 
S_{r+1}(q\cup\hatt{s})-S_{r+1}(\hatt{s}) =S_{r+1}(q)$
for every $q\in D_{\hatt{s}}$, a contradiction. 

Hence $r'\le r$, and the remaining claims are then obvious.
\end{proof}

\begin{proof}[Proof of Theorem~\ref{Thm:Extreme_symmetric_multiaffine_function_given_per_sums}]
We may assume $a<b$, put $n\coloneqq \#I$, and 
\[
 \widetilde{f}(p) &\coloneqq& f\big(a+(b-a)p_i:i\in I\big) \quad\text{ for }p\in[0,1]^I.
\]
Then $\widetilde{f}$ is a permutation invariant and separately affine-linear function on $[0,1]^I$,
so that 
Lem\-ma~\ref{Lem:Symmetric_partially_affine_functions} yields a function $g:\underline{n}_0\rightarrow\R$ 
with $\widetilde{f}(p)=\mathrm{BC}_p \,g$ for $p\in[0,1]^I$.
Hence, using also 
$
  S_\ell\big(a+(b-a)p_i:i\in I\big) = \sum_{k=0}^\ell\binom{\ell}{k}a^{\ell-k}(b-a)^kS_k(p)
$ for $\ell\in\underline{r}$ and $p\in[0,1]^I$,
and Lemma~\ref{Lem:S,K,E} with $\{A,B\}=\{S,K\}$, the claim follows from 
Theorem~\ref{Thm:Chebyshev-Hoeffding_k}.
\end{proof}

\section*{Acknowledgements}
We thank Norbert Henze and Patrick van Nerven for helpful remarks.

{\footnotesize
\section*{References}
Links to full texts, if provided here, are {\color{OliveGreen} green} if they are {\color{OliveGreen}free},
and {\color{BrickRed} red}  if  {\color{BrickRed} paywalled}, as experienced by us 
at various times 
starting in March 2022.
The asterisk *  marks 
references 
I have taken from secondary sources, without looking at the original.
Links in {\color{Blue} blue} are for navigating within this paper.

\renewcommand\refname{}

\end{document}